\documentclass[10pt]{amsart}

\title{The Crepant Resolution Conjecture}
\date{December 22nd, 2006}
\address{
Dept of Math, Univ. of British Columbia,
Vancouver, BC, Canada 
}
\email{jbryan@math.ubc.ca}
\address{
Dept of Math, California Institute of Technology, Pasadena, CA
}
\email{graber@caltech.edu}

\author{Jim Bryan and Tom Graber}

\usepackage{amsmath}
\usepackage{amsmath,amsthm,amsfonts}
\usepackage[all]{xypic}

\DeclareMathOperator{\Aut}{Aut}
\DeclareMathOperator{\Hilb}{Hilb}
\DeclareMathOperator{\Sym}{Sym}
\DeclareMathOperator{\GHilb}{G-Hilb}

\newcommand{\comment}[1]{}
\newcommand{\fz}{{\mathfrak{z}}}

\newcommand{\E}{{\mathbb E}}
\newcommand{\C}{{\mathbb C}}
\newcommand{\Q}{{\mathbb Q}}
\newcommand{\cnums} {{\mathbb C}}          
\newcommand{\znums} {{\mathbb Z}}       
\newcommand{\Z}{{\mathbb Z}}
\newcommand{\qnums} {{\mathbb Q}}       
\newcommand{\DD}{{\mathcal D}}
\newcommand{\bD}{{\bar{D}}}

\newcommand{\sym}{\operatorname{Sym}}
\newcommand{\hilb}{\operatorname{Hilb}}
\newcommand{\point}{\mathrm{pt}}

\newcommand{\cK}{{\mathcal{K}}}
\newcommand{\cC}{{\mathcal{C}}}
\newcommand{\oh}{{\mathcal{O}}}

\newcommand{\proj}{{\mathbb{P}}}
\renewcommand{\P}{\mathbb{P}}
\newcommand{\M}{\overline{{M}}}

\newcommand{\X}{\mathcal{X}}
\newcommand{\orb}{\mathit{orb}}

\newcommand{\vir}{{\mathrm{vir}}}
\newcommand{\IX}{{I\X}}

\newcommand{\NS}{\widehat{NS}}
\newcommand{\hbeta}{{\widehat{\beta}}}
\newcommand{\hF}{{\widehat{F}}}

\newcommand{\Hn}{\hilb^n(\C^2)}

\newcommand{\barH}{{\overline{H}}}
\newcommand{\Irr}{\operatorname{Irr}}
\newcommand{\Conj}{\operatorname{Conj}}
\newcommand{\age}{\operatorname{age}}

\newtheorem{thm}{Theorem}[section]

\newtheorem{cor}[thm]{Corollary}

\newtheorem{lemma}[thm]{Lemma}

\newtheorem{proposition}[thm]{Proposition}

\newtheorem{defn}[thm]{Definition}

\newtheorem{conjecture}[thm]{Conjecture}
\newtheorem{conjecture*}{Conjecture}

\newtheorem{rem1}[thm]{Remark}

\newenvironment{remark}{\begin{rem1}\em}{\end{rem1}}

\newcommand{\bpf}{\mathsc{Proof:}}
\newcommand{\epf}{\qed}

\begin{document}
\maketitle

\begin{abstract}
For orbifolds admitting a crepant resolution and satisfying a hard
Lefschetz condition, we formulate a conjectural equivalence between
the Gromov-Witten theories of the orbifold and the resolution.  We
prove the conjecture for the equivariant Gromov-Witten theories of
$\sym ^{n}\cnums ^{2}$ and $\hilb ^{n}\cnums ^{2}$.
\end{abstract}



\section{Introduction}

\subsection{Overview} Gromov-Witten theory is the mathematical
counterpart of topological string theory in physics. A well known
principle in physics states that string theory on an orbifold is
equivalent to string theory on a crepant resolution
\cite{Vafa-orbifold-numbers,Zaslow}. In their ground breaking paper
\cite{Chen-Ruan}, Chen and Ruan define orbifold cohomology using an
orbifold version of Gromov-Witten theory. Orbifold Gromov-Witten
theory was developed in the algebro-geometric context in
\cite{Ab-Gr-Vi-2,Ab-Vi,Ab-Gr-Vi} using Abramovich and Vistoli's notion
of twisted stable maps to a Deligne-Mumford stack.

The Gromov-Witten invariants of a projective manifold $Y$ are
multilinear functions $\left\langle \dotsb \right\rangle_{g,\beta
}^{Y}$ on its cohomology $H^{*} (Y)$. The Gromov-Witten invariants of
a orbifold $\X$ are multilinear functions $\left\langle \dotsb
\right\rangle_{g,\beta }^{\X }$ of the orbifold cohomology
$H^{*}_{\orb} (\X)$. Orbifold cohomology is by definition the ordinary
cohomology of the inertia stack $I\X $ with a shifted grading
\cite{Ab-Gr-Vi,Chen-Ruan}. $H_{\orb }^{*} (\X )$ contains the usual
cohomology $H^*(X)$ as a subspace and its orthogonal complement is
referred to as the space of twisted sectors.  If $\X $ is a Gorenstein
orbifold whose coarse moduli scheme $X$ admits a crepant resolution $Y
\to X$, Yasuda has proven that $H^{*}_{\orb} (\X ,\C)$ and $H^{*} (Y,
\C)$ are isomorphic as graded vector spaces \cite{Yasuda}.  Yasuda's
proof provides an equality of Betti numbers but does \emph{not}
provide any natural choice of isomorphism.  Nevertheless, it has been
suggested by Ruan \cite{Ruan-crepant} that there should be such an
isomorphism which identifies the Gromov-Witten theories.  He proposes
that specializing the values of certain quantum parameters of the
small quantum cohomology of the resolution will recover the orbifold
cohomology of the orbifold.  In this paper, we formulate an analogous
conjecture at the level of the full genus zero quantum potentials, and
we explore its consequences. In particular, we show this conjecture
allows one to essentially recover the Gromov-Witten theory of the
resolution in terms of that of the orbifold. We confirm the validity
of our conjecture for some examples including the case of $\X =\sym
^{n}\cnums ^{2}$, $Y=\hilb^{n}\cnums ^{2} $.

Recent work of Coates, Corti, Iritani, and Tseng, \cite{CCIT} strongly
suggests that for orbifolds failing the hard Lefschetz condition, the
relationship between the Gromov-Witten theories of the orbifold and
its resolution is more complicated and is better expressed in the more
sophisticated framework of Givental's Lagrangian formalism.  We are
grateful to them for bringing the hard Lefschetz condition to our
attention.


\subsection{Notation} By an \emph{orbifold}, we will mean a smooth
algebraic Deligne-Mumford stack over $\cnums $. An orbifold $\X $ is
said to be \emph{Gorenstein} if $\X $ has generically trivial
stabilizers and the canonical bundle of $\X $ pulls back from a line
bundle on the coarse moduli space $X$ (equivalently, for every $x\in
\X $, the action of the isotropy group on the canonical line bundle is
trivial).  A resolution of singularities $\pi :Y\to X$ is called
\emph{crepant} if $K_{Y}=\pi ^{*}K_{X}$.

Let $\X $ be a Gorenstein orbifold and let
\[
\pi :Y\to X
\]
be a crepant resolution of the coarse moduli space $X$.
\comment{For
simplicity of exposition, we will assume $\X $ and $Y$ have no odd
cohomology.} We say that an integer basis for the second homology group
of a variety is \emph{positive} if the cone generated by the basis
contains the Mori cone.

Let $\{\beta_1, \ldots , \beta_r \}$ be a positive basis of $H_2(Y)$
such that $\{\beta_{s+1}, \ldots ,\beta_r \}$ is a basis for the
kernel of $\pi _{*}:H_2(Y) \to H_2(X)$.  Note that $\{\pi _{*}\beta
_{1},\dotsc ,\pi _{*}\beta _{s} \}$ is a positive basis for $H_{2}
(X)$.  We choose a basis $\{\gamma_0, \ldots, \gamma_a\}$ for
$H^*(Y)$, with $\gamma_0=1$ and $\gamma_1, \ldots, \gamma_r$ the basis
for $H^2(Y)$ dual to the $\beta_i$.

The genus zero Gromov-Witten invariants of $Y$ are multilinear
functions $\left\langle \dotsb  \right\rangle_{\beta }^{Y}$ on $H^{*}
(Y)$, defined by cohomological evaluations against $[\M _{0,n}
(Y,\beta )]^{\vir }$, the virtual fundamental class of the moduli space
of stable maps \cite{coxkatz,mirrorbook}. The invariants are encoded
in the potential function,
\[
F^{Y} (y_{0},\dotsc ,y_{a},q_{1},\dotsc ,q_{r}) =\sum _{n_0,\ldots , n_a=0}^{\infty }\sum _{\beta } \left\langle
  \gamma_0^{n_0}\cdots\gamma_a^{n_a} \right\rangle_{\beta }^{Y}
\frac{y_0^{n_0}}{n_0!}\cdots\frac{y_a^{n_a}}{n_a!}q_1^{d_1 }\cdots
q_r^{d_r}
\]
where $\beta =d_{1}\beta _{1}+\dotsb +d_{r}\beta _{r}$ is summed over
all non-negative $(d_{1},\dotsc ,d_{r})$.

Similarly, the genus zero Gromov-Witten invariants of $\X $ are
multilinear functions $\left\langle \dotsb \right\rangle_{\beta }^{\X
}$ on $H^{*} _{\orb }(\X )$, defined by cohomological evaluations
against $[\M _{0,n} (\X ,\beta )]^{\vir }$, the virtual fundamental
class of the moduli space of \emph{twisted} stable maps
\cite{Ab-Gr-Vi}. We choose a basis $\{\delta _{0},\dotsc ,\delta _{a}
\}$ for $H^{*}_{\orb } (\X )$ and we define the potential function for
$\X $:
\[
F^{\X} (x_{0},\dotsc ,x_{a},u_{1},\dotsc ,u_{s}) =\sum
_{n_0,\ldots , n_a=0}^{\infty }\sum _{\beta } \left\langle
  \delta_0^{n_0}\cdots\delta_a^{n_a} \right\rangle_{\beta }^{\X}
\frac{x_0^{n_0}}{n_0!}\cdots\frac{x_a^{n_a}}{n_a!}u_1^{d_1 }\cdots
u_s^{d_s}
\]
where $\beta =d_{1}\pi _{*}\beta _{1}+\dotsb +d_{s}\pi _{*}\beta _{s}$
is summed over all non-negative $(d_{1},\dotsc ,d_{s})$.

The inertia stack $I\X $ of an orbifold $\X $ is defined to be the
fibered product of $\X $ with itself over the diagonal in $\X \times
\X $. The points of $I\X $ are pairs $(x,g)$ where $x\in \X $ and
$g\in \Aut _{\X } (x)$. There is an involution $I$ of $I\X $ taking
$(x,g)$ to $(x,g^{-1})$. To each component $\X _{i}$ of $I\X $ we
assign a rational number $\age (\X _{i})$ as follows. Let $(x,g)$ be a
point in $\X _{i}$. Then $g$ acts on $T_{x}\X $ with eigenvalues
$(\alpha _{1},\dotsc ,\alpha _{n})$ where $n=\dim \X $. Let $r$ be the
order of $g$ and define $s_{j}\in {0,\dotsc ,r-1}$ by $\alpha
_{j}=\exp (2\pi i\frac{s_{j}}{r})$. Then age is defined by
\[
\age (\X _{i}) = \frac{1}{r}\sum _{j=1}^{n}s_{j}.
\]
Age is well defined and is integral for Gorenstein orbifolds. 

As a graded vector space, the orbifold cohomology of $\X $ is the
cohomology of $I\X $ with the grading shifted by twice the age:
\[
H^{*}_{\orb } (\X ) = \bigoplus _{\X _{i}\subset I\X }H^{*+2\age (\X _{i})} (\X _{i}).
\]

Suppose that the coarse moduli space $X$ is projective with hyperplane
class $\omega $. In \cite{Fernandez}, Fernandez asked if the hard
Lefschetz isomorphism holds in orbifold cohomology, namely if the
operator $L_{\omega }$ given by multiplication by $\omega $ in the
\emph{orbifold cohomology ring}, induces isomorphisms
\[
L^{p}_{\omega }:H^{n-p}_{\orb } (\X )\to H^{n+p}_{\orb } (\X ).
\]
Fernandez proved that $L^{p}_{\omega }$ is an isomorphism for all
$\omega $ if and only if the age is invariant under the involution
$I$. We call this condition (also defined for non-projective
orbifolds) the hard Lefschetz condition.
\begin{defn}\label{defn: hard lefschetz}
An orbifold $\X $ is said to satisfy the hard Lefschetz condition if the involution 
\[
I:I\X \to I\X 
\]
preserves the age.
\end{defn}

Note that this condition is satisfied by holomorphic symplectic orbifolds.

\subsection{The Conjecture}

Our main conjecture relates the two
potential functions $F^{Y}$ and $F^{\X }$.

\begin{conjecture}[Crepant Resolution Conjecture]\label{mainconjecture}
Given an orbifold $\X$ satisfying the hard Lefschetz condition and
admitting a crepant resolution $Y$, there exists a graded linear
isomorphism
\[
L : H^*(Y) \to
  H_\orb^*(\X)
\]
and roots of unity $c_{s+1}, \ldots , c_r $ such that the
following conditions hold.
\begin{enumerate}
\item The inverse of $L$ extends the map $\pi ^{*}: H^*(\X)\to H^{*}
      (Y)$.
\item Regarding the potential function $F^Y$ as a power series in
      $y_{0},\dotsc ,y_{a},q_{1},\dotsc ,q_{s}$, the coefficients
      admit analytic continuations from $(q_{s+1},\dotsc ,q_{r}) =
      (0,\dotsc ,0)$ to $(q_{s+1},\dotsc ,q_{r}) = (c_{s+1},\dotsc
      ,c_{r})$.
\item The potential functions $F^\X$ and $F^{Y}$ are equal after
      the substitution
\begin{align*}
y_i &= \sum _{j} L^j_ix_j\\
q_i &= \begin{cases} c_i  & \text{when } i> s\\
u_i & \text{when } i\leq s. \end{cases}
\end{align*}
\end{enumerate}
\end{conjecture}

\begin{remark}\label{rem: FY and FX contain equiv info}
The cohomological parameters in the potential functions,
$\{x_{0},\dotsc x_{a} \}$ and $\{y_{0},\dotsc ,y_{a} \}$ are equal in
number by Yasuda's result. However, the number of quantum parameters,
$\{u_{1},\dotsc ,u_{s} \}$ and $\{q_{1},\dotsc ,q_{r} \}$ differ, and
so na\"ively, the potential function $F^{Y}$ appears to have more
information than $F^{\X }$. However, the divisor equation implies
that the potential function $F^{Y}$ contains redundant information. In
fact, given $L$, $c_{s+1},\dotsc ,c_{r}$, and $F^{\X }$, one can
essentially recover $F^{Y}$. This will be made more clear in
section~\ref{sec: degree in twisted sectors and QH}, where we will
present an alternative but equivalent formulation of the conjecture
which is particularly convenient when studying the small quantum
cohomology ring.  Rather than resolving the difference in the number
of $q$'s and $u$'s by setting some of the $q$'s to constants, one can
adjoin extra $u$ variables to the orbifold partition function, by
defining a generalized notion of degree for orbifold curves with
unmarked twisted points.
\end{remark}

\comment{\begin{remark}
The analytic continuation condition above is implicitly assumed in the
physical literature on the subject, and it is certainly reasonable to
believe that this even stronger convergence properties are true for
the potential functions in general, outside the context of crepant
resolutions of orbifolds.
\end{remark}}

\begin{remark}\label{rem: unstable terms}
A finite set of coefficients in these potential functions are not well
defined since certain degenerate moduli spaces do not exist. Namely,
terms of degree zero in the quantum parameters and of degree less than
three in the cohomological variables are undefined.  We are not making
any conjectures about these coefficients.  To get a precise equality,
one needs to either take triple derivatives of the series on both
sides, or choose compatible assignments of values to the unstable
coefficients.  It would be interesting to find a meaningful way of
defining these unstable invariants.
\end{remark}

\begin{remark}\label{rem: poincare pairing is preserved by L}
It is a consequence of the conjecture that the linear map $L$ must
preserve the (orbifold) Poincar\'e pairing.
\end{remark}

\begin{remark}\label{rem: can add descendants etc to conj}
If $\X$ admits an action of an algebraic torus $T$ and $Y$ is a
$T$-equivariant crepant resolution, then we can extend the conjecture
to include equivariant parameters. In fact, this equivariant 
version of the conjecture follows immediately from the absolute
version, by considering the conjecture applied to finite dimensional
approximations to the homotopy quotients $Y_{T}\to X_{T}$.
\end{remark}

\begin{remark}\label{rem: galois automorphisms}
The coefficients of $F^{Y}$ and $F^{\X }$ are rational
numbers, but in general the linear transformation $L$ may be have to
be defined over some extension of $\qnums $. A consequence of the
conjecture is that there is a symmetry of $F^{Y}$ given by the
action of the Galois group of the extension on the change of
variables. In practice, this is often a highly non-trivial symmetry.
\end{remark}

\begin{remark}\label{rem: CRC in higher genus}
Our conjecture may also hold as stated for higher genus
potentials. There is very little evidence in positive
genus, although Maulik's computation of the full Gromov-Witten potential for
$A_{n}$ surface resolution \cite{Maulik-An} does provide some positive
evidence.  The relationship of the higher genus Gromov-Witten
potentials for Gorenstein orbifolds failing the hard Lefschetz
condition is expected to be more complicated involving a mixing of
different genera. For a physical account, see
\cite{Aganagic-Bouchard-Klemm}. See also the discussion in section~5
of \cite{CCIT}.
\end{remark}

\subsection{The noncompact case}\label{subsec: noncompact case}

Although Gromov-Witten theory is best known in the compact setting,
the simple examples we want to focus on are non-compact, so we
observe that there is a large class of noncompact examples where
there is a well defined version of the conjecture.

The most convenient hypothesis here is to assume that $X$ is
projective over an affine scheme and that $Y$ is projective over $X$
and hence also projective over an affine.  (In fact, in our examples
$X$ will actually be affine.) In this setting, although the spaces
of stable (twisted) maps need not be proper, the evaluation maps
from the space of maps to $Y$ (or $I\X$) will be proper.  Thus we
have well-defined Gromov-Witten classes
\[
 \langle \gamma_1, \ldots , \gamma_n, *\rangle_\beta
\]
defined as in \cite{Ab-Gr-Vi-2} by pushing forward from the space of $n+1$
pointed (twisted) stable maps to $Y$ (or $I\X$).

If the target is projective, then because of the formula
\[
\langle
\gamma_1, \ldots , \gamma_{n+1}\rangle_{\beta} =
 \gamma_{n+1} (\langle \gamma_1, \ldots, \gamma_n, *\rangle_{\beta})
\]
these homology
valued invariants contain equivalent information to the numerical
Gromov-Witten invariants, and moreover, the conjecture as stated
implies immediately a conjecture for a generating series of
homology valued invariants.  In the noncompact setting, where one
cannot reduce the homology classes to numbers in this way, we can
then use these invariants to make a meaningful version of the
conjecture.

In fact, we will not need to pursue a careful language for these
refined invariants, since our examples have another useful feature.
They all admit a torus action with compact fixed locus.  Because of
this, there is a perfect pairing on the $T$-equivariant cohomology
given by formally applying the Bott residue formula.  While this
pairing takes values in $H_T^*(pt, \C)_{\mathfrak m}$, rather than
$\C$, it still allows us to do calculations at the level of the
familiar generating functions for numerical Gromov-Witten invariants
with the slight novelty that some of these numbers will be rational
functions in the equivariant parameters.  

\comment{A further advantage of this
equivariant case is that standard facts about equivariant cohomology
make it easy to conclude that if $\X$ admits an equivariant
compactification $\overline\X$ with equivariant crepant resolution
$\overline Y$ extending $Y$, then the conjecture for $\overline Y \to
\overline X$ implies this conjecture for $Y\to X$.  Since it is hard
to believe that Gromov-Witten theory would detect an obstruction to
the existence of such a resolution of singularities, it seems
reasonable to assume that if the conjecture holds for all projective
$X$, it will also hold for this nice class of quasi-projective
examples.}

\section{Degree in twisted sectors and Quantum cohomology}\label{sec: degree in twisted sectors and QH}

In this section we extend the definition of $\left\langle \dotsb
\right\rangle^{\X }_{g,\beta }$ to allow for $\beta $ to be a ``curve
class'' in the twisted sector. Consequently, the corresponding
Gromov-Witten potential of $\X $ includes quantum parameters
corresponding to twisted sectors. This allows us to formulate an
alternative version of the Crepant Resolution Conjecture where the
number of variables for $\X $ and for $Y$ are the same. In particular,
the large and small quantum cohomology rings of $Y$ and $\X $ have the
same number of deformation parameters and are isomorphic (in a certain
sense -- see subsection~\ref{subsec: orbifold quantum cohomology})
when the Crepant Resolution Conjecture holds.

\subsection{The orbifold Neron-Severi group and twisted degrees}

We define an enlarged Neron-Severi group for a Gorenstein orbifold
$\X$ as follows.  Let $T^1(\X)$ be the twisted part of
$H^2_{orb}(\X,\Z)$.  As this is generated by fundamental classes of
certain irreducible components of the inertia stack, it comes with a
canonical (unordered) basis and is a free Abelian group of rank $r-s$.

\begin{defn}
We define the \emph{orbifold Neron-Severi group}  $\NS _1(\X) $ by
\[
\NS _1(\X) = NS_1(\X) \oplus T^1(\X)^\vee.
\]
That is, an element $\hbeta \in\NS _1(\X)$ is a curve class $\beta $
in $\X$ together with a function $\hbeta (i)$ assigning an integer to
each age one component of $I\X$.  An
element of $\NS _1(\X)$ will be considered {\em effective} if the
underlying curve class is effective, and the function is nonnegative.
\end{defn}

Recall that evaluation at the $i$th point of a twisted stable map
takes values in the inertia stack and defines a virtual morphism
$e_{i}:\M _{g,n} (\X ,\beta )\to I\X $.  (We are using different
conventions here than those of \cite{Ab-Gr-Vi} or \cite{Ab-Gr-Vi-2}
--- our $\M$ corresponds to $\cK$ and our $e_{i}$ corresponds to
    $\tilde{e}_{i} $ of Proposition~6.1.4 of \cite{Ab-Gr-Vi-2}.)

\begin{defn}\label{defn: moduli space for twisted degrees}
Given an effective class $\hat \beta \in \NS_1(\X)$, we define
$\M_{g,n}(\X,\hat\beta)$ to be the moduli space parameterizing genus
$g$ twisted stable maps to $\X$ with degree $\beta$ with $n$ ordered
marked points and with $\hat\beta(i)$ \emph{unordered} twisted points which
map to $D_{i}$, the $i$th component of the inertia stack.  Precisely,
if we consider the following fiber product:
$$\xymatrix{
\M \ar[r] \ar[d] &
\overline{I}\X^n \times \overline{D}_1^{\hbeta(1)} \times \cdots \times \overline{D}_{r-s}^{\hbeta(r-s)} \ar[d]\\
\M_{g,n+\sum \hat\beta(i)}(\X,\beta)\ar[r] & \overline{I}\X^n \times
\overline{I}\X^{\hbeta(1)} \times \cdots \times
\overline{I}\X^{\hbeta(r-s)}}$$ then we define $\M_{g,n}(\X,\hbeta)$
to be the quotient $[\M/S_{\hbeta(1)} \times \cdots \times
S_{\hbeta(r-s)}].$ Here $\overline{I}\X $ is the rigidified stack and
$\overline{D}_{i}$ is the $i$th component.
\end{defn}

\begin{remark}\label{rem: counting interpretation of degree in twisted sector}
One interpretation of the usual degree is as counting the number of
times a curve intersects some fixed divisor. Similarly, we can
interpret the degree in the twisted sector as counting the number of
times some curve ``intersects'' some twisted divisor, namely it gives
the number of (unmarked, non-nodal) stacky points that get mapped to
the corresponding age one component of the inertia stack. The reason
for not including nodal stacky points in the count is so that the
degree will be locally constant in families.  We ignore the marked
points so that degree is additive when gluing smooth curves together
to form nodal ones, and therefore the boundary of the moduli spaces
$\M _{g,n} (\X ,\hbeta )$ have a product description analogous to the
usual one for the ordinary stable map moduli spaces.
\end{remark}

\subsection{Gromov-Witten invariants for degrees in twisted sectors
and the divisor equation}

Since the right hand vertical arrow of the diagram in
Definition~\ref{defn: moduli space for twisted degrees} is simply an
inclusion of a union of connected components, so is the left hand
vertical arrow, which means that the perfect obstruction theory and
virtual fundamental class for the usual space of twisted stable maps
immediately give one on $\M$, and by descent, we get a virtual
fundamental class on $\M_{g,n}(\X,\hbeta)$.

We can thus define Gromov-Witten invariants for curves with degrees
defined in the twisted sectors using these moduli
spaces. Correspondingly, we define the genus zero \emph{extended
Gromov-Witten potential} of $\X $ by
\begin{align*}
\hF^{\X } (x_{0},\dotsc ,x_{a},u_{1},\dotsc ,u_{r}) &=\\
\sum
_{n_0,\ldots , n_a=0}^{\infty }\sum _{\hbeta  }& \left\langle
  \delta_0^{n_0}\cdots\delta_a^{n_k} \right\rangle_{\hbeta }^{\X}
\frac{x_0^{n_0}}{n_0!}\cdots\frac{x_a^{n_a}}{n_a!}u_1^{d_1 }\cdots
u_s^{d_s}u_{s+1}^{\hbeta  (1)}\dotsb u_{r}^{\hbeta (r-s)}
\end{align*}

The extended invariants do not contain any new information, since we
have the following obvious formula:
\begin{equation}\label{eqn: div eqn for twisted divisors}
\langle D_1^{\hbeta(1)} \cdots D_{r-s}^{\hbeta(r-s)}\alpha _1 \cdots
\alpha _{n} \rangle_{\beta }^{\X } = \hbeta(1)!\dotsb \hbeta (r-s)! \cdot
\langle \alpha_1 \cdots \alpha_n\rangle_{\hbeta }^{\X }
\end{equation}
which immediately reduces the calculation of these ``new'' invariants
to the calculation of the standard orbifold invariants.
We think of this as the analog of the divisor equation for the
``twisted divisors'' $D_i$, since it formally allows us to remove
the $D_i$ from invariants.  Note, however, that this equation
is different in form from the usual divisor equation.

It is useful to see what the divisor equation tells us about the form
of the potential function.  For $Y$, it is well known that repeated
application of the divisor equation implies that (up to unstable
terms) we have

\comment{
\[
F^Y =  \sum _{n_{r+1},\ldots , n_a=0}^{\infty
}\sum _{\beta } \left\langle
  \gamma_{r+1}^{n_{r+1}}\cdots\gamma_a^{n_k} \right\rangle_{g,\beta }^{Y}
\frac{y_{r+1}^{n_{r+1}}}{n_{r+1}!}\cdots\frac{y_a^{n_a}}{n_a!}
\left(q_{1}e^{y_{1}} \right)^{d_{1}}\dotsb \left(q_{r}e^{y_{r}} \right)^{d_{r}}
\]}

\[
F^Y=F^Y(y_0,0,0,\ldots, 0,y_{r+1},\ldots,y_{a},q_1e^{y_1},\ldots , q_re^{y_r}).
\]

In other words, the potential function depends on the variables in the
combinations
\[
q_1e^{y_1},\dotsc ,
q_re^{y_r}, y_{r+1},\dotsc , y_n.
\]
We can apply this to only the exceptional classes, giving the form more useful to us here:
\[
F^Y=F^Y(y_0,\ldots, y_s,0,\ldots, 0,y_{r+1},\ldots,y_{a},q_1,\ldots ,
q_s,q_{s+1}e^{y_{s+1}},\ldots, q_re^{y_r})
\]
For the orbifold invariants, the analogous result is that the extended
potential function depends on the variables only in the combinations
\[
u_1e^{x_1},\dotsc , u_se^{x_s}, (u_{s+1} + x_{s+1}),\dotsc ,
(u_r+x_r), x_{r+1},\dotsc ,x_n.
\]
More precisely, equation~\eqref{eqn: div eqn for twisted divisors}
implies the identity

\[
\hF^\X= F^\X(x_0 , \ldots , x_s, (x_{s+1}+u_{s+1}),
\ldots , (x_r+u_r), x_{r+1}, \ldots , x_n,u_{1},\dotsc ,u_{s}).
\]

So, assuming Conjecture \ref{mainconjecture}, we see that we get the
equality $\hF^\X = F^Y$ for the extended potential function after the
change of variables:

\begin{align*}
y_i &= \sum _{j} L^j_ix_j\\
q_i &= \begin{cases} c_ie^{L_i^ju_j}  & \text{when } i> s\\
u_i & \text{when } i\leq s. \end{cases}
\end{align*}

Since this change of variables is invertible up to the discrete
choices of branches of certain logarithms, it shows that one can
essentially recover the Gromov-Witten theory of $Y$ from that of $\X$.
Moreover, in this form it is especially clear that the existence of
the standard divisor equation on $Y$ gives a very strong and
mysterious prediction about the potential for $\X$ -- it should depend
on the new $u$ variables (or equivalently some of the original $x$
variables) only in terms of certain exponentials.

\subsection{Orbifold Quantum Cohomology}\label{subsec: orbifold
quantum cohomology} Another application of the above formalism is to
define a quantum product for an orbifold that is equivalent to the
quantum product of its crepant resolution by a method completely
parallel to the usual definition.  We will discuss here only the small
quantum cohomology.  Of course, one can use the derivatives of the
genus zero potential function to define a big quantum cohomology ring
for orbifolds and everything we say can be applied there as well.

Assume $\X$ is a Gorenstein orbifold with projective coarse moduli scheme.
We consider the three evaluation maps from $\M_{0,3}(\X,\hbeta)$ to
$\IX$, and given classes $\delta$ and $\gamma$ in $H^*(\IX)$, we
define
\[
\delta * \gamma= \sum_\hbeta \left(\langle \delta, \gamma, *\rangle_{\hbeta }^{\X }
\right)^{\vee }u^\hbeta
\]
where $(\cdot )^{\vee }$ denotes dual with respect to the orbifold
Poincar\'e pairing.  The same proof of
associativity holds for this product as for the one considered in
\cite{Ab-Gr-Vi}.

We can express the quantum product in a basis, using the orbifold
Poincar\'e pairing $ g_{ij}$ on $H^{*} (\IX )$ as
\[
\gamma * \delta = \sum \langle \gamma, \delta, \gamma_i
\rangle_{\hbeta}\, g^{ij}\gamma_j u^\hbeta.
\]

Hence, it is an immediate consequence of
Conjecture~\ref{mainconjecture} that the products agree in the sense
that if we identify $H^*(Y)$ and $H^*(I\X)$ using $L$, then the
structure constants for the quantum product are related by the change
of variables:
\begin{align*}
q_i &\mapsto \begin{cases} c_ie^{L_i^ju_j}  & \text{when } i> s\\
u_i & \text{when } i\leq s. \end{cases}
\end{align*}

\begin{remark}As in subsection~\ref{subsec: noncompact case}, this
definition of the quantum product makes sense using only the
hypothesis that $X$ is projective over an affine scheme.  The argument
reducing the equivalence of the quantum products of $\X$ and $Y$ to
the equivalence of the potential functions uses the perfectness of the
Poincar\'e pairing, which we do have in the torus equivariant
setting provided that the fixed locus is compact.

\end{remark}


\begin{remark}
The definition of small quantum cohomology given in \cite{Chen-Ruan}
or \cite{Ab-Gr-Vi} can be recovered from this one by setting the new
parameters equal to zero.  It follows then, that one recovers that
quantum cohomology ring of $\X$ from the quantum cohomology ring of
$Y$ by simply setting some of the $q$'s to roots of unity. The idea of
setting quantum parameters on the resolution equal to roots of unity
first appears in the mathematics literature in the work of Ruan
\cite{Ruan-crepant} where he observes that in some examples, one needs
to set $q=-1$ to recover the orbifold cohomology of $\X $.
\end{remark}

\section{Examples}

To provide evidence for our conjecture we consider orbifolds of the
form
\[
\X =\left[V/G \right]
\]
where $G\subset SL (V)$ is a finite subgroup.

When the dimension of $V$ is 2 or 3, there is a canonical crepant
resolution given by the $G$-Hilbert scheme \cite{BKR}:
\[
Y=\GHilb (V).
\]
The diagonal $\cnums ^{\times}$ action on $V$ commutes with $G$ and the
induced action on $X$ lifts to $Y$. Thus the crepant resolution
conjecture can be considered $\cnums ^{\times}$ equivariantly.

By \cite{BKR}, there is a canonical basis for $H^{*}_{\cnums ^{\times}}
(\GHilb V) $ indexed by $R\in \Irr (G)$, irreducible representations
of $G $. On the other hand, there is a canonical basis of
$H^{*}_{\cnums ^{\times},\orb } ([V/G])$ indexed by $(g)\in \Conj (G)$,
conjugacy classes of $G$. Denote the corresponding cohomology variables by 
\[
\left\{y_{R} \right\}_{R\in \Irr (G)} \text{ and }\left\{x_{(g)}
\right\}_{(g)\in \Conj (G)}
\]
respectively. Let $y_{0}$ and $x_{0}$ be the variables corresponding
to the trivial representation and the trivial conjugacy class
respectively.

\subsection{Polyhedral and Binary polyhedral groups.}

A finite subgroup $G$ of $ SO (3)$ (respectively $SU (2)$) is called a
\emph{polyhedral} (respectively \emph{binary polyhedral}) group. Such
groups are classified by ADE Dynkin diagrams and they come with a
natural representation $V$ of dimension 3 (respectively 2). For these
groups, the equivariant quantum cohomology of $\GHilb (V)$ has been
completely described in terms of the root theory of the corresponding
ADE root system by Bryan-Gholampour
\cite{Bryan-Gholampour2,Bryan-Gholampour3}. They conjecture that the
change of variables for the crepant resolution is a certain
modification of the character table.

\begin{conjecture}\label{conj: change of vars for polyhedral G}
The change of variables for the crepant resolution conjecture in the
case of 
\[
\GHilb (V) \to V/G
\]
where $G$ is a polyhedral or binary
polyhedral group is given by
\begin{align*}
y_{0}&=x_{0},\\
y_{R}& = \frac{1}{|G|}\sum _{g\in G} \sqrt{\chi _{V} (g)-\dim V}\quad
\chi _{R} (g) \, x_{(g)},\\
q_{R}&=\operatorname{exp}\left(\frac{2\pi i \dim R}{|G|} \right)
\end{align*}
where $R$ runs over the \emph{non-trivial } irreducible
representations of $G$.
\end{conjecture}
Note that as a consequence of $V$ being the natural representation of
a polyhedral or binary polyhedral group, the orbifold $\X =[V/G]$
satisfies the hard Lefschetz condition. Moreover, all non-trivial
conjugacy classes have age one, and so the above linear transformation
preserves the grading. In fact, these are the only faithful group
representations that have the property that all non-trivial elements
have age one.

Using the root theoretic formula for the Gromov-Witten potential of
$\GHilb (V)$ given in \cite{Bryan-Gholampour2,Bryan-Gholampour3} and
applying the crepant resolution to the above change of variables, one
arrives at a prediction for the orbifold Gromov-Witten potential
$F^{\X }$. This prediction has not been verified in general, but it
does pass some strong tests of its validity. Namely, it can be shown
to exhibit various vanishing properties and to have the correct
classical terms.

The complete determination of $F^{\X }$, and hence the verification of
the crepant resolution conjecture, has been done for $G$ equal to
\begin{align*}
\znums _{2}&\subset SU (2) \text{ in the next subsection,}\\
\znums _{3}&\subset SU (2) \text{ in \cite{Bryan-Graber-Pandharipande},}\\
\znums _{4}&\subset SU (2) \text{ in \cite{Bryan-Jiang},}\\
\znums _{2}\times \znums _{2}&\subset SO (3) \text{ in \cite{Bryan-Gholampour1}, and}\\
A_{4}&\subset SO (3) \text{ in \cite{Bryan-Gholampour1}.}
\end{align*}

\subsection{The case of the rational double point.}\label{subsec:
C2modZ2 case} We consider the case where $V=\cnums ^{2}$ and $G=\{ \pm
1\}\subset SU (2)$ so that
\[
\X =[\cnums ^{2}/\{\pm 1\}],\quad  Y=T^{*}\P ^{1}.
\]
This is the simplest nontrivial example and it already provides a very
interesting case study. Here we will establish the equivariant version
of the conjecture.

Let $T=\C^{\times } \times \C^{\times }$ so that
\[
H^{*}_{T} (\point )\cong \qnums [t_{1},t_{2}].
\]
The natural $T$ action on $\cnums ^{2}$ induces a $T$ action on $Y$,
the minimal resolution of the quotient $X=\cnums ^{2}/\{\pm 1 \}$. $Y$
is isomorphic to $T^{*}\P ^{1}$, the total space of the cotangent
bundle of $\P ^{1}$. There are two fixed points of the $T$ action on
$Y$ having weights
\[
(2t_{1},t_{2}-t_{1})\text{ and } (2t_{2},t_{1}-t_{2}).
\]

First we will compute the genus zero potential function for $Y$.  We
take our generator for $H_2(Y)$ to be the class of the zero section,
$[E]$.  We let $\gamma \in H^2_T(Y)$ be the dual of $[E]$. It is given
by the first Chern class of an equivariant line bundle with weights
$-t_{1}$ and $-t_{2}$ at the fixed points. 

The degree zero invariants are simply given by triple intersections in
equivariant cohomology
\[
\left\langle a,b,c \right\rangle_{0} = \int _{Y}a\cup b\cup c
\]
which are computed by localization. The results are:

\[
\left\langle1,1,1\right\rangle_{0}= \frac{1}{2t_1t_2},\quad 
\left\langle\gamma, 1, 1\right\rangle_{0} = 0,\quad 
\left\langle\gamma, \gamma, 1\right\rangle_{0} = -\frac12,\quad 
\left\langle\gamma,\gamma,\gamma\right\rangle_{0}=\frac12 (t_1+t_2)
\]

To compute the invariants in positive degrees, we first observe
that the image of any nonconstant morphism from a curve to $Y$
must lie in $E$.  Thus we have a natural isomorphism
\[
\M_{0,n}(Y,d[E]) \cong \M_{0,n}(\proj^1, d), 
\]
however, the virtual fundamental classes on the two sides differ.
Under the above identification, it is well known that
\[
[\M_{0,n}(Y,d[E])]^\vir = e(R^1\pi_*f^*N_{E/Y})
\]
where $\pi :\mathcal{C}\to \M _{0,n} (\P ^{1},d)$ and
$f:\mathcal{C}\to Y$ are the universal curve and the universal map
respectively.

Since $E$ is a -2 curve, we have an isomorphism of
$N_{E/Y} \cong \oh(-2)$.  Consider the standard
Euler sequence on $\proj^1$,
\[
 0 \to \oh(-2) \to \oh(-1)\oplus\oh(-1) \to \oh \to 0 .
\]
Pulling this sequence back to $\cC$ and taking the associated long
exact sequence of derived pushforwards gives us
$$
0\to \oh \to R^1\pi_*f^*(\oh(-2)) \to R^1\pi_*f^*(\oh(-1)\oplus\oh(-1)) \to 0.
$$
An analysis of the weights shows that the action of $T$ on the left
hand term in this sequence is given by $t_{1}+t_{2}$.  We conclude
that
\[
e (R^1\pi_*f^*N_{E/Y})=(t_1+t_2)e(R^1\pi_*(\oh(-1)+\oh(-1))).
\]
The integral is then
evaluated using the famous Aspinwall-Morrison formula:
\[
\left\langle\,  \right\rangle_{d}= (t_{1}+t_{2})\int _{[\M _{0,0} (\P
^{1},d)]}e (R^{1}\pi _{*}f^{*} (\oh (-1)\oplus \oh
(-1)))=\frac{t_{1}+t_{2}}{d^{3}}.
\]

Let $y_{0}$ and $y_{1}$ denote the variables corresponding to $1$ and
$\gamma $. Combining the above formulas with the divisor equation and
the point axiom, we have shown the following.
\begin{proposition}\label{prop: FY for TP1}
The genus zero Gromov-Witten potential function of $Y$ is given by:
\[
F^Y = \frac{1}{12t_1t_2} y_0^3  - \frac14 y_0y_1^2 +
\frac{t_1+t_2}{12}y_1^3+ (t_{1}+t_{2})\sum_{d>0} \frac{1}{d^3} q^d e^{dy_1}.
\]
\end{proposition}

We now consider the invariants for the orbifold $\X$. Let $1$ and $D$
be generators for $H^{0}_{\orb } (\X )$ and $H^{2}_{\orb } (\X )$ and
let $x_{0}$ and $x_{1}$ be the corresponding variables.

Since the coarse
moduli space for $\X$ is affine, every stable map is constant.  If the
source curve has any twisted points, the image of the map is forced to
be the unique point of $\X$ with nontrivial stabilizer.  Thus we see
that with the exception of 
\[
\langle 1, 1, 1 \rangle=\frac{1}{2t_1t_2},
\]
every invariant naturally arises as an integral over $\M_{0,n}(B\znums
_2)$. By the point axiom and monodromy considerations, the only other
invariant involving 1 is
\[
\left\langle 1,D,D \right\rangle=\frac{1}{2}.
\]

Since the only remaining non-vanishing invariants are then $\langle
D^n \rangle$ we actually need only consider the connected component of
$\M _{0,n} (B\znums _{2})$ where all the evaluation maps go to the
twisted sector.  Setting $n=2g+2$, we denote this space as
$\barH_g^{ord}$. Concretely, it is the usual compactified moduli space
of hyperelliptic curves (with ordered branch points). The virtual
class on $\barH _{g}^{ord}$ is given by 
\[
e (R^{1}\pi _{*}f^{*} (L\oplus L))
\]
where $f$ and $\pi $ are the universal map and universal curve for $\M
_{0,n} (B\znums _{2})$ and
\[
L\oplus L\to B\znums _{2}
\]
is two copies of the non-trivial line bundle over $B\znums _{2}$ with
the torus acting with weight $t_{1}$ on the first factor and with
weight $t_{2}$ on the second factor. The bundle $R^{1}\pi _{*}f^{*}L$
is in fact isomorphic to $\E ^{\vee }$, the dual of the Hodge bundle
pulled back by the map $\barH _{g}^{ord}\to \M _{g}$. We conclude that
that for $n=2g+2>0$, we can write
\begin{align*}
\left\langle D^{n} \right\rangle &=\int _{\barH
_{g}^{ord}}e (\E ^{\vee }\oplus \E ^{\vee })\\
&=  - (t_{1}+t_{2})\int _{\barH _{g}^{ord}} \lambda _{g}\lambda _{g-1}.
\end{align*}

The generating function for these integrals was computed in
\cite[Corollary 2]{Faber-Pandharipande-logarithmic}. Applying that
computation, we obtain:
\begin{proposition}\label{prop: F for X=C2modZ2}
The potential function of $\X =[\cnums ^{2}/\{\pm 1 \}]$ is given by
\[
F^{\X } (x_{0},x_{1}) =
\frac{1}{12t_{1}t_{2}}x_{0}^{3}+\frac{1}{4}x_{0}x_{1}^{2} -
(t_{1}+t_{2})x_{1}^{2}H (x_{1})
\]
where (following the notation of
\cite{Faber-Pandharipande-logarithmic}) $H (x_{1})$ satisfies
\[
(x_{1}^{2}H (x_{1}))''' = \frac{1}{2}\tan \left(\frac{x_{1}}{2} \right).
\]
\end{proposition}
\begin{cor} The crepant resolution conjecture holds for the pair
$(Y,\X )$. That is, the potential functions $F^{Y} (y_{0},y_{1},q)$
and $F^{\X } (x_{0},x_{1})$ agree, up to unstable terms, under the
change of variables (c.f. Conjecture~\ref{conj: change of vars for
polyhedral G})
\[
y_{0}=x_{0},\quad y_{1}=ix_{1},\quad q=-1.
\] 
\end{cor}
\textsc{Proof:} Clearly the terms of $F^{Y}$ and $F^{\X } $ which have
$y_{0}$ and $x_{0}$ match up. And since we are only interested in
stable terms, it suffices to check that
\[
\left(\frac{d}{dx_{1}} \right)^{3}F^{Y} (x_{0},ix_{1},-1) =
\left(\frac{d}{dx_{1}} \right)^{3}F^{\X } (x_{0},x_{1}).
\]
The right hand side is given by
\[
- (t_{1}+t_{2})\frac{1}{2}\tan \left(\frac{x_{1}}{2} \right),
\]
whereas the left hand side is
\begin{align*}
&(t_{1}+t_{2})\left[\frac{i^{3}}{2}+\sum _{d=1}^{\infty }i^{3} (-e^{ix_{1}})^{d} \right]\\
=&(t_{1}+t_{2})\frac{1}{2i}\left[\frac{1-e^{ix_{1}}}{1+e^{ix_{1}}} \right]\\
=&(t_{1}+t_{2})\frac{1}{2}\tan \left(\frac{-x_{1}}{2} \right).
\end{align*}
\qed

\subsection{The case of the Hilbert scheme}\label{subsec: the case of 
Hilb}
We consider the case where
\[
\X =\Sym^{n} (\cnums ^{2})\text{ and  } Y=\Hilb^{n} (\cnums ^{2}).
\]
This is one of the best known and most studied examples of a crepant
resolution of singularities of a Gorenstein orbifold.

We will show in this section that by matching the Nakajima basis for
the cohomology of the Hilbert scheme with the natural basis for the orbifold
cohomology of the symmetric product we verify
Conjecture~\ref{mainconjecture} in this case.

Because the Hilbert scheme is holomorphically symplectic, there are no
interesting Gromov-Witten invariants unless one works equivariantly.
An analogous fact is true on the orbifold side.  Thus, to verify the
conjecture for the nonequivariant theory it suffices to compare the
ring structure on the ordinary cohomology of the Hilbert scheme with
the orbifold cohomology of the symmetric product. This is done in
\cite{Vasserot,Lehn-Sorger} (see also \cite{Fantechi-Gottsche,
Uribe}).

The nontrivial, fully equivariant genus 0 Gromov-Witten theory of $Y$ is
determined in
\cite{Ok-Pan-Hilb}.  We will determine the genus 0 equivariant
Gromov-Witten theory of $\X$ and verify that it matches their result after
the appropriate change of variables.

First, let us describe the inertia stack $I\X$.  We use the standard
correspondence between conjugacy classes of $S_n$ and partitions of
$n$.  Given such a partition $\mu$, the corresponding component of the
inertia stack $I_\mu$ can be described by choosing a representative
permutation $\sigma$ and taking the stack quotient $[\C^{2n}_\sigma /
C(\sigma)]$ where $\C^{2n}_\sigma$ denotes the invariant part of
$\C^{2n}$ under the action of $\sigma$ and $C(\sigma)$ denotes the
centralizer of $\sigma$ in $S_n$.  The dimension of $\C^{2n}_\sigma$
is $2l(\mu)$ and the age of $\mu$ is $n-l(\mu)$.  The quotient of a
vector space by a finite group has no higher cohomology groups, so we
conclude that a basis for the orbifold cohomology of $\X$ as a
$\Q[t_1,t_2]$-module is given by 
\[
[I_\mu]\subset H^{2n-2l(\mu)}_{T,\orb } (\X )
\]

Since each element of $S_n$ is conjugate to its inverse, the
equivariant Poincar\'e pairing on $H^*_{T,\orb}(\X)$ is diagonal in
this basis.  It is easily computed by localization, since the fixed
points of the $T$ action on $I\X$ are isolated --- there is a single
fixed point in each irreducible component $I_\mu$.  This point has
automorphism group equal to the centralizer of a representative
element, which has order
\[
\fz (\mu) = |\Aut (\mu)|
\prod \mu_i.
\]
It follows that the pairing is given by
\[
( [I_\mu],[I_\mu])=\frac{1}{\fz(\mu)}(t_1t_2)^{-l(\mu)}.
\]

It is straightforward to check that the orbifold product here is a slight
modification of the usual multiplication on $Z\qnums [S_{n}]$, the
center of the group ring of $S_n$ obtained by inserting factors of
$t_1t_2$ to make that product respect the grading by age.  In
particular, the limit $t_1=t_2=1$ gives the standard product on
$Z\qnums [S_{n}]$.

On the Hilbert scheme, there is an analogous description of the equivariant
cohomology, given by the Nakajima basis.  Given a partition $\mu$, the
corresponding class 
\[
N_\mu \in H^{2n-2l (\mu )}_T(\Hn)
\]
is given
by $\frac{1}{\prod \mu_i}[C_\mu]$ where $C_\mu$ is the
subvariety of $\Hn$ whose general point parameterizes a length $n$
subscheme composed of $l(\mu)$ irreducible components
of lengths $\mu_i$.  The $T$-equivariant Poincar\'e pairing
in the Nakajima basis is also diagonal with
\[
( N_\mu , N_\mu) = \frac{(-1)^{n-l(\mu)}}{\fz(\mu)}(t_1t_2)^{-l(\mu)}.
\]

This gives us an obvious candidate for the map $L$ identifying
the orbifold cohomology of $\Sym^n(\C^2)$ with the cohomology
of $\Hn$. Namely, we define $L$ by
\begin{equation}\label{eqn: change of vars for HilbC2}
L([I_\mu])= i^{l(\mu)}N_\mu.
\end{equation}
Note also, that since there exists a unique partition of length $n-1$,
the partition corresponding to a 2-cycle, which we will denote $(2)$,
there is only one divisor class, and so a single constant $c$ to
choose to finish determining the change of variables.  The correct
choice of $c$ turns out to be $-1$. Thus the predicted change of
variables for the quantum parameters is
\[
q=-e^{iu}.
\]

To establish the full equality of the genus zero Gromov-Witten
potentials, it will be extremely convenient to use the formalism
introduced for the small quantum product as a bookkeeping device.  Let
$c^{\nu }_{\mu }$ be the structure constants for quantum
multiplication by $[I_{(2)}]$:
\[
 [I_{(2)}]*[I_\mu] = \sum_\nu c^\nu_\mu [I_\nu]. 
\]
Here the $c_\mu^\nu$ are elements of $\Q[t_1,t_2][[u]]$ where $u$ is
the quantum parameter associated to the twisted sector as defined in
Section~\ref{sec: degree in twisted sectors and QH}.  If we let
$c_\mu^\nu(d)$ denote the coefficient of $u^d$ in $c^\nu_\mu$, then we
have the formula
\begin{equation}\label{gwforc}c_\mu^\nu(d)=\fz(\nu)(t_1t_2)^{l(\nu)}\langle
[I_\mu], [I_\nu], [I_{(2)}]\rangle_d.
\end{equation}
Note that the above Gromov-Witten invariant is an element of
$\Q(t_1,t_2)$, whereas $c^{\nu }_{\mu } (d)$ is a polynomial. This
fact will be essential for the degree arguments that follow.  While
the polynomality is an immediate consequence of the existence of the
equivariant quantum product referred to in Section~\ref{subsec:
noncompact case}, the reader can also check that it follows directly
from the explicit localization formula we will give in the next
section.

By degree considerations, we see that $c_\mu^\nu$ vanishes if
$l(\mu)\geq l(\nu)+1$.  Since Equation \ref{gwforc} gives a symmetry,
we also have the inequality $l(\nu)\geq l(\mu) +1$.  We will see below
that if $|l(\nu) - l(\mu)|=1$ then the only contribution to
$c_\mu^\nu$ is in degree zero where we just see the classical term
corresponding to multiplication in the group ring of $S_n$.  Aside
from these classical terms, the matrix $c^\nu_\mu$ is diagonal.

\begin{lemma}\label{diagonal} If $l(\mu)=l(\nu)$, but $\mu \neq \nu$,
then $c_\mu^\nu =0$.
\end{lemma}
\textsc{Proof:} Let $\M _{0,(\lambda ,\mu ,\nu )} (\X ,d)$ denote the
component(s) of $\M _{0,3} (\X ,d)$ given by $e_{1}^{-1} (I_{\lambda
})\cap e_{2}^{-1} (I_{\mu })\cap e_{3}^{-1} (I_{\nu })$.

By definition, we have
\[
c_\mu^\nu (d)[I_{\nu }] =
e_{3*}([\M_{0,((2),\mu,\nu)}(\X,d)]^\vir)^{\vee }.
\]
By degree consideration, this must be a codimension 1
class in $I_\nu$.  However it is easy to see that the codimension of
$e_{3}(\M_{0,((2),\mu,\nu)}(\X,d))$ is at least 2, since the
intersection of the images of $I_\mu$ and $I_\nu$ in $\X$ has
codimension at least 2 in each.  The lemma follows immediately.  \qed

To finish the determination of the structure of the quantum cohomology
ring, we use a localization calculation.  Because $\Sym^n(\C^2)$ is
affine, every twisted stable map is constant at the level of coarse
moduli schemes.  It follows that we have a canonical identification of
the $T$ fixed locus of the space of maps to $\X$ with the space of
maps to $BS_n$ (the fixed locus of the action of $T$ on $\X$).  The
normal bundle to this fixed locus decomposes naturally as a sum of two
rank $n$ vector bundles. These two bundles come with $T$ weights
$t_1$ and $t_2$, but are otherwise identical, each corresponding to
the standard $n$-dimensional representation of $S_n$ under the usual
correspondence between sheaves on $BG$ and representations of $G$.  We
will use $V$ to denote this bundle. There is another way to think of
$V$ which is convenient for us here. Consider the morphism
$i:BS_{n-1} \to S_n$ induced by the standard inclusion of $S_{n-1}
\hookrightarrow S_n$.  Then $V$ is simply the pushforward of the
structure sheaf.

Since the coarse moduli scheme of $BS_n$ is a point, the moduli
space of twisted maps $\M_{0,r}(BS_n)$ is smooth of dimension $r-3$.
By the results of \cite{Gr-Pa} we can identify the equivariant
virtual fundamental class of $\M_{0,r}(\Sym^n(\C^2))$ with the pushforward
from this fixed locus of the class
\[
e(-R^{\bullet } \pi_*f^*(V\oplus V)) .
\]
where the torus acts on the two factors of $V$ are with weights $t_1$ and
$t_2$, and where $f$ and $\pi$ are the universal maps in the universal
diagram
$$\xymatrix{\mathcal{C} \ar[r]^f \ar[d]^\pi & BS_n \\
\M_{0,r}(BS_n).}
$$
To give a description of the virtual class in more familiar terms, we extend
the above diagram to

$$\xymatrix{ \tilde{\mathcal{C}} \ar@/_1pc/[dd]_p \ar[r]^{\tilde{f}} \ar[d]^g & BS_{n-1} \ar[d]^i \\
\mathcal{C} \ar[r]^f \ar[d]^\pi & BS_n \\
\M_{0,r}(BS_n)}.
$$

The curve $\tilde{\mathcal{C}}$ here is the degree $n$ covering of
$\mathcal{C}$ corresponding to the map to $BS_n$ via the usual correspondence
between principal $S_n$ bundles and degree $n$ \'etale covers.

We know that $V=i_*\oh_{BS_{n-1}}$.  Since $i$ is a finite morphism,
we have $f^*(i_*\oh)=g_*\tilde{f}^*\oh_{BS_{n-1}}=g_*\oh_{\tilde{C}}$.
Since $g_*$ is exact, we conclude that the virtual class can be
rewritten as the pushforward from the fixed locus of
\[
e(-R^{\bullet } p_* (
\oh_{\tilde{\mathcal{C}}} \oplus \oh_{\tilde{\mathcal{C}}})).
\]
In other words, if we think of the space of maps to $BS_n$ as
parameterizing the family of $n$-sheeted covers of $\proj^1$ given by
$\tilde{\mathcal{C}}$, then the invariants we want to compute are
expressed in terms of the Chern classes of the Hodge bundle.  Thus if
we let
\[
\E^\vee=R^1p_*\oh_{\tilde{\mathcal{C}}}
\]
and we let $s$ denote the locally constant function on
$\M_{0,r}(BS_n)$ recording the number of connected components of the
fibers of $\tilde{\mathcal{C}}$, we obtain the following formula.

\begin{lemma}\label{locform} The $r$ point, degree zero invariants of
$\X $ are given by
\[
\langle \mu_1 , \ldots, \mu_r\rangle_{0} =\int_{[\M_{0,(\mu_1, \ldots
, \mu_r)}(BS_n)]^\vir} (t_1t_2)^{-s}c_{top}(\E_{t_1}^\vee \oplus
\E_{t_2}^\vee)
\]
where $\M _{0,(\mu _{1}\dotsb \mu _{r})} (BS_{n})$ is the component(s) of
$\M _{0,r} (BS_{n})$ where $e_{i}$ maps to the component of $IBS_{n}$
corresponding to $\mu _{i}$.  We also must interpret the above
integral as a sum over connected components of the moduli space. The
rank of $\E$ and the integer $s$ can vary from component to component.
\end{lemma}

Since the orbifold $\X $ only has divisor classes in the twisted
sector, the higher degree invariants are determined by the degree zero
invariants by the divisor equation. Thus the above lemma determines
all the invariants.

\begin{lemma} If $d>0$, then $c_\mu^\nu(d)$ is divisible by $(t_1+t_2)$.
\end{lemma}
This is a consequence of Mumford's relation that $c(\E \oplus \E^\vee)
= 1$.  If we set $t_2=-t_1$ and use the fact that
$c_{top}(\E^\vee_{-t_1}) = \pm c_{top}(\E_{t_1})$, we find that the
integrand in Lemma \ref{locform} is simply a power of $t_1$.  The
hypothesis $d>0$ implies that the moduli space is positive
dimensional, so the result follows. \qed

We remark that on the Hilbert scheme side, this divisibility is
related to the existence of a holomorphic symplectic structure on
$\Hilb^n(\C^2)$.

\begin{cor}If $d>0$ and $\mu \neq \nu$, we have $c^\nu_\mu=0$.
\end{cor}
Given Lemma \ref{diagonal} and the discussion just before it, we see
that the only interesting case here is if $l(\mu) = l(\nu)+1$.
However, in this case, we know that the degree of $c_\mu^\nu$ is zero,
so the divisibility constraint forces this invariant to vanish.\qed

This reduces our task to the calculation of the invariants
$c^\mu_\mu$.  We can further reduce to the case where the partition
$\mu$ has just one part, by using the following lemma.

\begin{lemma} 
$$c^\mu_\mu =\frac{1}{\fz(\mu)} \sum_{i=1}^{l(\mu)} \mu_i
c_{(\mu_i)}^{(\mu_i)}.$$
\end{lemma}

The right hand side of this formula is easily seen to be the
contribution from those components of $\M_{0,((2),\mu,\mu)}(\X,d)$
where the corresponding branched cover $C$ consists of $l(\mu)$
connected components, all but one of which is a smooth genus zero
curve branched only at 0 and $\infty$.  To prove the lemma, we need to
show that the other components make no contribution.  We will do this
by means of the formula of Lemma \ref{locform}, so we will always be
considering the space of maps to $BS_n$ rather than to $\X$, and we
will denote the two distinguished points of the source curve
corresponding to $\mu$ as 0 and $\infty$.  These are the only points
over which the associated branched cover of $\proj^1$ has non-simple
branching.

\smallskip

\noindent {\bf Step 1}: Suppose we have a component where the
associated cover is connected.  Then, by the Riemann-Hurwitz formula,
it will have genus $g=\frac{d+3}{2} - l(\mu)$.  In order for the
integral in Lemma \ref{locform} not to vanish, it is obviously
necessary that $2g\geq d$ since $d$ is the dimension of the moduli
space.  This inequality is satisfied only if $l(\mu)=1$ (in which case
the Lemma is vacuously true).  Otherwise, we conclude that a component
of the moduli space can contribute to this invariant only if it
parameterizes disconnected covers.

\smallskip

\noindent {\bf Step 2}: Suppose we consider a component $\M'$ of the
moduli space where the corresponding branched cover is disconnected.
We get a natural map $\Psi: \M' \to (\prod_a \M_a)/\Aut$ where the
$\M_a$ are some moduli spaces of lower degree branched covers with
certain branching conditions and the group $\Aut$ is acting by
permuting factors with identical parameters.  We do not need a very
careful description here, since we will use just two crude facts.
First, if two different factors of the target space parameterize covers
with branching away from zero and $\infty$, then $\Psi$ has positive
dimensional fibers, since we can independently act by $\C^{\times }$ on
different components.  Since $\E$ is pulled back under $\Psi$ this
immediately kills contributions from any such component of $\M$.

If a branched cover has all the simple branch points on a single
connected component, then the other components are necessarily genus
zero curves ramified only at zero and infinity.  Now Step 1 will apply
to the remaining interesting component, showing that this component of
moduli space makes no contribution to the integral unless we are in
the maximally disconnected case. \qed

We remark that the argument in this lemma extends to give an alternate
proof of Lemma \ref{diagonal}.

We see then, that we will have completely determined the quantum
multiplication by $[I_{(2)}]$ once we calculate the invariants
$c_{(n)}^{(n)}$ for all $n$.  Here we can give an explicit formula.

\begin{lemma} We have the following:
\[
 \sum_d \langle [I_{(n)}], [I_{(n)}], [I_{(2)}] \rangle_d u^d =
-\frac{i}{2}\frac{t_1+t_2}{t_1t_2}(n\cot(\frac{nu}{2}) - \cot(\frac{u}{2})).
\]
\end{lemma}

\textsc{Proof:} This follows from the same argument as \cite[Theorem
6.5]{Br-Pa-local-curves}.  The restriction of $\tilde{\mathcal{C}}$ to
$e_{1}^{-1} (I_{(n)})\subset \M_{0,3} (\X ,d)$ is necessarily a family
of connected curves of genus $g$ where $d=2g-1$.  Applying
Lemma~\ref{locform} and the divisor equation, we get the formula
\[
\langle [I_{(n)}], [I_{(n)}], [I_{(2)}] \rangle_d =
-\frac{t_1+t_2}{t_1t_2}\frac{1}{(2g-1)!}\int_{[\M_{0,((n)(n)(2)\dotsb (2))}
(BS_{n})]}\lambda_g\lambda_{g-1}.
\]
The map to $\M_{g,2}$ induced by the family $\tilde{\mathcal{C}}$ is
generically finite of degree $(2g)!$ onto its image, which is the set
of curves admitting a degree $n$ map to $\proj^1$ totally ramified at
the two marked points. The image of this map is called
$\overline{H}_{d}\subset \M _{g,2}$ in \cite{Br-Pa-local-curves} and
the pairing of $[\overline{H}_{g}]$ against $\lambda_g\lambda_{g-1}$
is explicitly evaluated in \cite{Br-Pa-local-curves} to yield the
series above.  \epf

Having completely determined the $c^\nu_\mu$ we can deduce our main result.
\begin{thm} \label{thm: QH(Sym)=QH(Hilb))}
After making the change of variables given by
equation~(\ref{eqn: change of vars for HilbC2}) and relating the
quantum parameters by 
\[
q=-e^{iu},
\] 
the genus zero Gromov-Witten potential of $\Hilb^n(\C^2)$ is equal to
the (extended) genus zero Gromov-Witten potential of
$\Sym^n(\C^2)$. Hence the crepant resolution conjecture holds in this
case.
\end{thm}

By direct inspection, the matrix of multiplication by $[I_{(2)}]$ in
$QH^*(\Sym^n(\C^2))$ matches with the matrix of multiplication by
$i[N_{(2)}]$ in $QH^*(\Hilb^n(\C^2))$ calculated in \cite{Ok-Pan-Hilb}
(equations (6) and (8), see also \cite{Br-Pa-local-curves} equations
(19) and (29)) under the change of variables $q=-e^{iu}$.  As is
observed there, the fact that this matrix has distinct eigenvalues
implies that after extending the scalars to $\Q(t_1,t_2)$ the quantum
cohomology is generated by the divisor class $[I_{(2)}]$.  Thus the
entire ring structure is encoded in this multiplication matrix.

Finally, since the small quantum cohomology is generated by divisors,
a variant of the reconstruction theorem of Kontsevich-Manin shows that
one can use the WDVV equation to reduce arbitrary genus zero
Gromov-Witten invariants to invariants with only two insertions
(c.f. \cite{Rose}). As these are already encoded in the small quantum
product, the proof of the theorem is complete. \qed

\pagebreak

\subsection{Equivalence with other theories}

There are two other theories which are equivalent to the quantum
cohomologies of $\sym ^{n} (\cnums ^{2})$ and $\Hilb ^{n} (\cnums
^{2})$. By computing the equivariant Gromov-Witten partition function
(in all genus) for the degree $n$ invariants of $\P ^{1}\times \cnums
^{2}$ relative to $\{0,1,\infty \}\times \cnums ^{2}$, one obtains the
structure constants of an associative Frobenius algebra
\cite{Br-Pa-local-curves}. Similarly, on obtains a Frobenius algebra
from the partition function for the degree $n$ equivariant
Donaldson-Thomas invariants of $\P ^{1}\times \cnums ^{2}$ relative to
$\{0,1,\infty \}\times \cnums ^{2}$
\cite{Okounkov-Pandharipande-dtlc}. Theorem~\ref{thm:
QH(Sym)=QH(Hilb))} completes the following tetrahedron of
equivalences.

\vspace{1.5in}
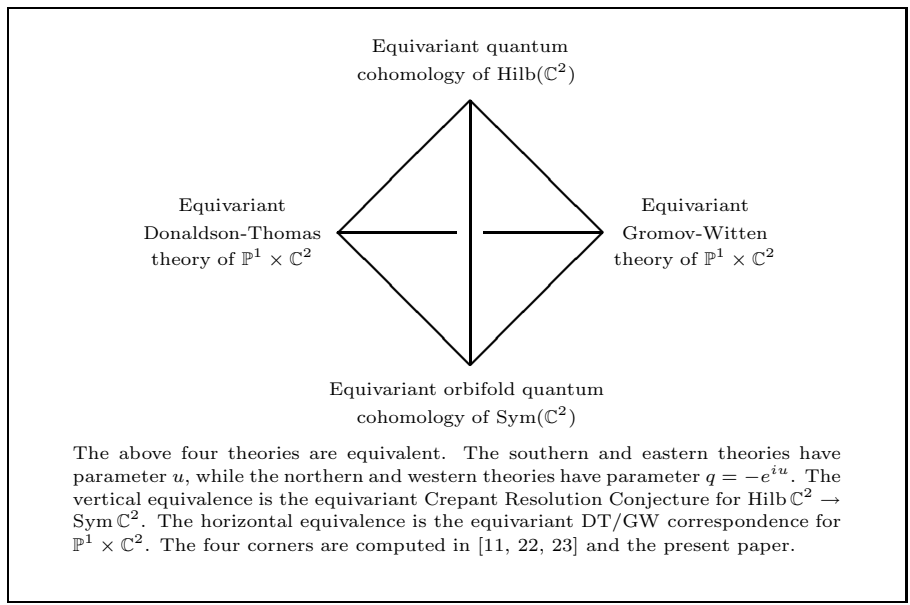
\begin{center}
\scriptsize
\begin{picture}(200,75)(-30,-50)
\put(-100 ,-115 ){\line(1 ,0 ){340}}
\put(-100 ,-115 ){\line(0 ,1 ){225}}
\put(-100 ,110 ){\line(1 ,0 ){340}}
\put(240 ,-115 ){\line(0 ,1 ){225}}
\thicklines
\put(25,25){\line(1,1){50}}
\put(25,25){\line(1,-1){50}}
\put(125,25){\line(-1,1){50}}
\put(125,25){\line(-1,-1){50}}
\put(75,-25){\line(0,1){100}}
\put(25,25){\line(1,0){45}}
\put(80,25){\line(1,0){45}}
\put(75,95){\makebox(0,0){Equivariant quantum}}
\put(75,85){\makebox(0,0){cohomology of $\Hilb (\cnums ^{2})$ }}
\put(75,-35){\makebox(0,0){Equivariant orbifold quantum }}
\put(75,-45){\makebox(0,0){cohomology of $\operatorname{Sym} (\cnums ^{2})$ }}
\put(160,35){\makebox(0,0){Equivariant}}
\put(160,25){\makebox(0,0){Gromov-Witten}}
\put(160,15){\makebox(0,0){theory of $\P ^{1}\times \cnums ^{2}$}}
\put(-15,35){\makebox(0,0){Equivariant}}
\put(-15,25){\makebox(0,0){Donaldson-Thomas}}
\put(-15,15){\makebox(0,0){theory of $\P ^{1}\times \cnums ^{2}$}}
\end{picture}
\end{center}
\begin{quote}
\scriptsize{The above four theories are equivalent. The southern and
eastern theories have parameter $u$, while the northern and western
theories have parameter $q=-e^{iu}$.  The vertical equivalence is the
equivariant Crepant Resolution Conjecture for $\Hilb \cnums ^{2}\to
\operatorname{Sym}\cnums ^{2}$. The horizontal equivalence is the
equivariant DT/GW correspondence for $\P ^{1} \times \cnums ^{2}$. The
four corners are computed in \cite{Br-Pa-local-curves,
Okounkov-Pandharipande-dtlc,Ok-Pan-Hilb} and the present paper.  }
\end{quote}
\vspace{.5in}

\subsection{Acknowledgments.}  We warmly acknowledge helpful
discussions with Mina Aganagic, Renzo Cavalieri, Tom Coates, Amin Gholampour,
Yunfeng Jiang, Rahul Pandharipande, Michael Thaddeus, and Hsian-Hua
Tseng. We acknowledge support from NSERC, NSF, the Sloan Foundation,
and IHES.

\bibliography{mainbiblio}
\bibliographystyle{plain}

\end{document}